\documentclass[11pt]{article}

\usepackage{amsmath,amssymb,amsthm,mathtools,geometry}
\usepackage{booktabs}
\usepackage{tabularx}
\usepackage{array}
\usepackage{enumitem}
\usepackage{verbatim}
\usepackage{xcolor}
\usepackage{etoolbox}
 
\usepackage[
    colorlinks=true,
    linkcolor=blue,
    citecolor=green,
    urlcolor=blue
]{hyperref}
\newcolumntype{Y}{>{\raggedright\arraybackslash}X}

\newtheorem{theorem}{Theorem}[section]
\newtheorem{lemma}{Lemma}[section]
\newtheorem{proposition}{Proposition}[section]
\newtheorem{corollary}{Corollary}[section]
\newtheorem{assumption}{Assumption}[section]

\theoremstyle{definition}
\newtheorem{definition}{Definition}[section]

\theoremstyle{remark}
\newtheorem{remark}{Remark}[section]

\newcommand{\R}{\mathbb R}
\newcommand{\C}{\mathbb C}

\newcommand{\calS}{\mathcal S}
\newcommand{\calW}{\mathcal W}
\newcommand{\aff}{\operatorname{aff}}

\renewcommand{\d}{{\rm d}}
\numberwithin{equation}{section}
\newcommand{\revisioncolor}{red}

\newcommand{\setreferencecolor}[1]{%
  \color{black}%
  \ifstrequal{#1}{UnnikrishnanHaghighatshoarVetterli:2018}
    {\color{\revisioncolor}}{}%
  \ifstrequal{#1}{TsakirisPeng:2019}
    {\color{\revisioncolor}}{}%
  \ifstrequal{#1}{PengWangTsakiris:2021}
    {\color{\revisioncolor}}{}%
}

\title{On the Unique Identification of Finitely Many Radial Heat Sources\thanks{B. Jin is partly supported by Hong Kong RGC General Research Fund (14306824) and ANR / Hong Kong RGC Joint Research Scheme (A-CUHK402/24), and NSFC / RGC Joint Research Scheme (N\_CUHK446/25). The work of Z. Zhou is supported by Hong Kong Research Grants Council (15303122) and an internal grant of Hong Kong Polytechnic University (Project ID: P0038888, Work Programme: 1-ZVX3).}}

\author{Bangti Jin\thanks{Department of Mathematics, The Chinese University of Hong Kong, Shatin, New Territories, Hong Kong, P.R. China (\texttt{bangti.jin@gmail.com, b.jin@cuhk.edu.hk}).}   \and Zhi Zhou\thanks{Department of Applied Mathematics, The Hong Kong Polytechnic University, Kowloon, Hong Kong, P.R. China (\texttt{zhizhou@polyu.edu.hk})}}
\date{}

\begin{document}

\maketitle

\begin{abstract}
In this work, we investigate an inverse source problem for the heat equation in $\R^d$ ($d\geq1$), in which the source is a finite positive linear combination of translates of a known real-valued radial Schwartz profile with nonzero total mass. Both the the locations and strengths of the sources are to be identified from the observations of time traces at finitely many spatial points over an arbitrary nonempty open positive-time interval. For $d\geq2$, we prove that the sharp universal measurement threshold is $dN+1$: every set of $dN+1$ measurement points in affine general position uniquely identifies every source configuration with at most $N$ locations. For $d=1$, we prove that two distinct measurement points uniquely identify every finite source configuration, independently of the number of sources, and that this count is sharp. The analysis relies on an equivalence between equality of heat traces and equality of weighted distance measures, and then employs perpendicular-bisector hyperplanes to resolve geometric ambiguities. We also establish identification guarantees using $d+1$ measurement points: clustered designs under a reference-distance separation condition, adaptive designs for arbitrary configurations, and any fixed affinely independent design for almost every location configuration, uniformly over all positive strengths.\\
\textbf{Key words}: inverse source problem, radial source, uniqueness, universal sharp threshold
\end{abstract}

\section{Introduction}\label{sec:intro}

Inverse source problems are concerned with recovering an unknown source in an evolution equation from observational data, and arise in heat conduction, diffusion, contaminant transport, and pollution source localization. These problems typically lack uniqueness unless one imposes suitable structural conditions on the source, e.g., separable sources or point sources \cite{CannonPerez:1991,ElBadiaHaDuong:2002}. In this work, we consider a geometrically natural class of sources consisting of finitely many translates of a common radial profile. Let $d\geq1$, and let $\rho\in\calS(\R^d)$ be a known real-valued radial Schwartz function with nonzero total mass. For a finite set $A\subset\R^d$ and positive strengths $q=(q_a)_{a\in A}$, the source takes the form
\begin{equation*}
f_{A,q}(x)=\sum_{a\in A}q_a\rho(x-a),
\end{equation*}
where both the source locations in $A$ and their strengths are unknown. This model describes finitely many localized sources with a common spatial shape and possibly different strengths. 
Approximate point sources can also be represented by normalized and rescaled profiles concentrating at the source locations. Given the configuration $f_{A,q}$, let $u=u_{A,q}$ solve
\begin{equation}\label{eq:heat}
\left\{
\begin{aligned}
\partial_tu-\Delta u &=f_{A,q},
&&\text{in }\R^d\times(0,\infty),\\
u(\cdot,0)&=0,
&&\text{in }\R^d.
\end{aligned}
\right.
\end{equation}
The observations are given by $u_{A,q}(b,t)$ for $b\in B$ and $t\in J$, where $B\subset\R^d$ is a finite measurement set and $J\subset(0,\infty)$ is a fixed nonempty open interval. Thus each measurement location provides a time trace. The main question is: \emph{how many spatial measurement points guarantee unique identification of the locations and strengths of every source configuration with at most $N$ locations, for every measurement set satisfying a natural nondegeneracy condition?} The nondegeneracy condition refers to affine general position (AGP): every subset of at most $d+1$ measurement points is affinely independent; see Definition \ref{def:AGP}.

In this work, we establish both universal guarantees over all AGP measurement designs and reduced sufficient bounds under additional design or genericity conditions. In Theorem \ref{thm:main}, we prove that, for $d\geq2$, the sharp universal threshold is $dN+1$: every AGP design of size $dN+1$ uniquely identifies the locations and positive strengths of every configuration with at most $N$ sources, whereas some AGP design of size $dN$ fails even for unit-strength configurations. Importantly, identifiability holds among all finite configurations, with no cardinality bound imposed on the competing source set. For $d=1$, in Theorem \ref{thm:main-1d}, we  show that two distinct measurement points suffice to identify every finite configuration, and that this count is sharp. We also obtain substantially reduced sufficient bounds: $d+1$ clustered measurements identify configurations with separated squared distances from a reference point;  $d+1$ adaptively chosen measurements identify every finite configuration; and any fixed AGP design of $d+1$ points identifies almost every location configuration of each finite cardinality, for every choice of positive strengths. These results are given in Theorems \ref{thm:clustered}, \ref{thm:adaptive}, and \ref{thm:genericdesign}, respectively. See Section \ref{sec:main} for the precise statements and further discussions.

The proof of the universal bounds is divided into an analytic part and a geometric part. First, we establish an exact analytic--geometric equivalence: equality of heat traces at a measurement point $b$ on any nonempty open positive-time interval is equivalent to equality of the weighted distance measures $\nu_b(A,q)=\sum_{a\in A}q_a\delta_{|b-a|}$. Thus each heat trace determines the observed distances and the total strength at each distance, but does not distinguish sources equidistant from the sensor. The analytic implication follows from holomorphic continuation of a Gaussian transform, recovery of centered radial moments, and a triangular moment argument followed by inversion of a finite Vandermonde system. Geometrically, positivity of the source strengths and equality of the weighted distance measures ensure that every competing source location must match a distance occurring in the original configuration at each sensor. We employ a covering argument using perpendicular-bisector hyperplanes and AGP to identify the source locations, and then determine the strengths using measurements separating individual distances. We prove sharpness using a unit-strength reflection construction based on two interlaced regular polygons and a generic selection lemma for points on finitely many hyperplanes. In one dimension, composing reflections about two measurement points rules out any nonzero finitely supported ambiguity. Finally, the reduced sufficient bounds resolve the ambiguity across sensors in three different ways: clustered measurements preserve the ordering of separated reference distances; adaptive measurements first restrict all possible source locations to a finite set and then separate them with one additional sensor; and fixed AGP designs exclude inconsistent distance assignments outside algebraic sets of measure zero. The exceptional sets in the last argument depend only on the source locations and the measurement design, not on the positive strengths.

Inverse source problems for parabolic equations have been studied extensively, with classical treatments addressing uniqueness, stability, and reconstruction from boundary, final-time, or interior observations; see the monographs \cite{Cannon:1984,Isakov:1990}. For point-source models, uniqueness has been investigated for heat and advection--dispersion--reaction equations
\cite{ElBadiaHaDuong:2002,ElBadiaHaDuongHamdi:2005},
including the identification of multiple moving sources from boundary
observations \cite{AndrleElBadia:2012}.  Uniqueness from spatially sparse observations has also been studied. Hettlich and Rundell \cite{HettlichRundell:2001} studied the identification of discontinuous spatial sources in the heat equation on the unit disc from flux measurements at two carefully chosen boundary points; see also \cite{RundellZhang:2020,LiZhang:2020,LinZhangZhang:2022} for further developments. More recently, the works \cite{GuZhangZhang:2025,GongJinKianLiu:2026,DalmassoCenKian:2026} investigated the recovery of a point source or a more general measure-valued source for the heat equation in unit balls in $\mathbb{R}^d$ from sparse boundary measurements at suitably chosen points. See also \cite{KusiakWeatherwax:2007,BenBelgacem:2012} for the identification of a moving point source from sparse measurements. The present work differs from these existing studies in both the source model and the measurement setting. We study finite positive linear combinations of translates of a known radial profile and measure pointwise time traces at finitely many fixed locations in the whole space, rather than boundary fluxes in a bounded domain. This distinction also changes the structure of the inverse problem: boundary observations reflect the geometry of the domain, whereas the whole-space model encodes source--sensor distances. The works \cite{GuZhangZhang:2025, GongJinKianLiu:2026, DalmassoCenKian:2026} rely crucially on the spectral structure of the Dirichlet Laplacian on the unit disc or ball, including angular eigenfunctions and Bessel functions. In contrast, we exploit an equivalence between heat-trace data and weighted distance measures and then uses Euclidean geometry to derive identification guarantees and sharp universal measurement thresholds. The approach also applies to finite sums of Dirac point sources with positive strengths, provided that no measurement point coincides with a source location, for which the analytic reduction to weighted distance measures is more direct, while the geometric arguments remain unchanged.

In two dimensions, Ling et al. \cite{LingYamamotoHonTakeuchi:2006} studied a closely related setting with unit source strengths. They derived distance-moment identities and stated a uniqueness result for three noncollinear measurement points \cite[Theorem 1]{LingYamamotoHonTakeuchi:2006}. In the radial unit-strength model, however, the trace at each sensor determines only an unlabeled multiset of source distances. Establishing uniqueness still requires resolving the consistency of distance assignments across sensors. Proposition \ref{prop:4pointcounterexample} gives an explicit ambiguity in this setting, which shows that noncollinearity alone is insufficient for uniqueness. The present work revisits this issue by treating the unlabeled distance data explicitly in the more challenging setting of arbitrary unknown positive source strengths. Theorems \ref{thm:main} and \ref{thm:main-1d} provide sharp measurement thresholds for uniform identification over all AGP designs in every dimension $d\geq1$. Theorems \ref{thm:clustered}, \ref{thm:adaptive}, and \ref{thm:genericdesign} provide further identification guarantees with reduced sufficient bound under additional assumptions on the measurement design or the source configuration.

The geometric reduction connects our problem to Euclidean distance geometry; see \cite{LibertiLavorMaculan:2014} for a survey. It is also related to reconstruction from unlabeled or partially labeled distance data. Boutin and Kemper \cite{BoutinKemper:2004} showed that generic finite point configurations can be recovered, under suitable hypotheses, up to rigid motion and relabeling from the distribution of their pairwise distances, while exceptional ambiguities can occur. Our data differ from the pairwise-distance distribution of an unknown configuration. Here, the measurement points have known absolute positions, and each heat trace determines a weighted distance measure from one of these known points to the unknown sources. Thus the data are labeled by sensor but not by source. The central geometric issue is to determine whether these separate weighted distance measures admit more than one common source configuration.

The rest of the paper is organized as follows. In Section \ref{sec:main}, we present the main results and discuss their connections and distinctions. In Section \ref{sec:direct}, we review the direct problem and derives the moment identities. In Section \ref{sec:equivalence}, we establish the analytic--geometric equivalence underlying the analysis. In Section \ref{sec:bounds}, we prove the universal identification guarantees and their sharpness through geometric upper and lower bounds. Finally, in Section \ref{sec:design}, we establish the clustered-design, adaptive, and almost-everywhere identification results.

\section{Main results and discussions}\label{sec:main}

Throughout, we impose the following condition on the profile $\rho$, and all source strengths are positive.

\begin{assumption}\label{ass:profile}
The function $\rho\in\calS(\R^d)$ is real valued and radial, with
\[
M_0:=\int_{\R^d}\rho(x)\,\d x\neq0.
\]
\end{assumption}

The measurement points in $B$ satisfy the following condition.

\begin{definition}\label{def:AGP}
A finite set $B\subset\R^d$ is in \textit{affine general position} (AGP) if every subset of $B$ of cardinality at most $d+1$ is affinely independent.
\end{definition}

\subsection{Universal guarantees}
In this part, we present several universal guarantees and will give the detailed proofs in Section \ref{sec:bounds}.
The first result gives the sharp measurement threshold uniform
over all AGP designs for $d\geq 2$.
\begin{theorem}\label{thm:main}
Let Assumption \ref{ass:profile} hold, and let $d\geq2$ and $N\geq1$.  The smallest integer $m$ such that every AGP set of $m$ measurement points uniquely identifies every configuration with at most $N$ sources among all finite configurations is $m=dN+1$. More precisely:
\begin{enumerate}
\item[{\rm(i)}]
Let $B\subset\R^d$ be in AGP and satisfy $|B|\geq dN+1$.
Let $(A,q)$ and $(A',q')$ be configurations, with $|A|\leq N$.
If \(u_{A,q}(b,t)=u_{A',q'}(b,t)\) for every \(b\in B\) and \(t\in J\),
then \(A=A'\) and \(q_a=q'_a\) for every \(a\in A\).

\item[{\rm(ii)}]
There exist distinct $N$-point sets $A,A'\subset\R^d$ and an AGP set
$B\subset\R^d$, with $|B|=dN$, such that
\(u_{A,\mathbf1}(b,t)=u_{A',\mathbf1}(b,t)\) for every \(b\in B\)
and \(t>0\), where $\mathbf1$ denotes unit strengths.
\end{enumerate}
Consequently, $dN+1$ is the smallest measurement count for which every
AGP design identifies every configuration with at most $N$ sources
among all finite configurations.
\end{theorem}

The number $dN+1$ is a worst-case threshold over all AGP
measurement sets. The lower bound follows by constructing one specially
arranged set of $dN$ measurements that fails. Theorem \ref{thm:main} does not rule out the possibility that a generic or deliberately
designed set of $dN$, or even fewer, measurement points identifies all or most source configurations.
It is worth noting that the cardinality bound applies only to $A$, with no bound imposed on the competing configuration $A'$. The proof combines the exact equivalence between heat traces
and weighted distance measures in
Section \ref{sec:equivalence} with positivity of the strengths
and geometric counting arguments based on perpendicular-bisector
hyperplanes in Section \ref{sec:bounds}.

The second result gives a sharp measurement count when $d = 1$.

\begin{theorem}\label{thm:main-1d}
Let Assumption \ref{ass:profile} hold, let $d=1$, and let $b_1,b_2\in\R$ be distinct. Let $(A,q)$ and $(A',q')$ be configurations. If $u_{A,q}(b_j,t)=u_{A',q'}(b_j,t)$ for $j=1,2$ and every $t\in J$, then $A=A'$ and $q_a=q'_a$ for every $a\in A$. 
\end{theorem}

Theorem \ref{thm:main-1d} shows that when $d=1$, for every $N\geq1$,
two distinct measurement points identify every configuration
with at most $N$ source locations among all finite
configurations. Unlike the threshold $dN+1$ for $d\geq2$,
this sharp count is independent of $N$.

For a single source in any dimension, Theorems \ref{thm:main} and \ref{thm:main-1d} give the following sharp bound.

\begin{corollary}\label{cor:singleton}
Let Assumption \ref{ass:profile} hold, and let $d\geq1$. Any $d+1$ affinely independent measurement points uniquely identify every single-source configuration among all finite configurations. No choice of fewer than $d+1$ fixed measurement points guarantees this identification for every single-source configuration.
\end{corollary}

\subsection{Design-dependent and adaptive measurement bounds}

The threshold $dN+1$ in Theorem \ref{thm:main} is a worst-case guarantee over every AGP measurement
set and every source configuration with at most $N$ locations. In practice, substantially fewer measurements
suffice if additional geometric information is available, if measurements can be chosen adaptively,
or if an exceptional set of configurations is permitted. These cases are covered by the
following three results. The proofs of these results are provided in Section \ref{sec:design}.

The first result gives the unique identification for a fixed design with \(d+1\) clustered measurements
under the hypothesis of separation of the squared source distances from a reference point. Then
the measurements sufficiently close to that point preserve the distance ordering and resolve the labeling ambiguity.
Specifically, for $b_0\in\R^d$, $R>0$, and $\delta>0$, we define the set
\begin{equation*}
\begin{aligned}
\mathfrak A(R,\delta;b_0)
=\bigl\{
A\subset\overline{B(b_0,R)}:\;
A\text{ is finite, and }
\bigl||a-b_0|^2-|\widetilde a-b_0|^2\bigr|\geq\delta
\text{ for all distinct }a,\widetilde a\in A
\bigr\}.
\end{aligned}
\end{equation*}
The separation condition is vacuous when $|A|\leq1$.

\begin{theorem}\label{thm:clustered}
Let Assumption \ref{ass:profile} hold, and let $d\geq1$, $b_0\in\R^d$, $R>0$, and $\delta>0$. Choose $b_1,\ldots,b_d\in\R^d$ such that $b_0,\ldots,b_d$ are affinely independent and $\max_{1\leq j\leq d}|b_j-b_0|<\frac{\delta}{4R}$. Let $(A,q)$ and $(A',q')$ be configurations, with $A\in\mathfrak A(R,\delta;b_0)$. If $u_{A,q}(b_j,t)=u_{A',q'}(b_j,t)$ for $j=0,\ldots,d$ and every $t\in J$, then $A=A'$ and $q_a=q'_a$ for every $a\in A$.
\end{theorem}

Note that in Theorem \ref{thm:clustered}, the radius and separation conditions apply only to $A$, not to the competing location set $A'$. The measurement count $d+1$ is sharp among fixed designs, since the prescribed class contains every single-source configuration in $B(b_0,R)$ and the reflection argument in Corollary \ref{cor:singleton} applies. Although the sensor count does not explicitly depend on $|A|$, the class parameters imply $|A|\leq\lfloor \frac{R^2}{\delta}\rfloor+1$. When the reference distances are distinct, the parameters needed to choose the remaining sensors can be recovered from the first trace; see Corollary \ref{cor:adaptiveclustered}. Further, we give an explicit reconstruction formula for consistently labeled distances in Remark~\ref{rem:labelledreconstruction}.

The next result requires neither a separation condition nor a cardinality bound.

\begin{theorem}\label{thm:adaptive}
Let Assumption~\ref{ass:profile} hold, and let $d\geq1$.
Every finite configuration in $\R^d$ can be uniquely identified among all finite configurations using $d+1$ adaptively chosen measurements. The first $d$ measurement points may be any affinely independent points. Based on their traces, the $(d+1)$st point can be chosen so that the combined observations uniquely identify the configuration and the full measurement set is affinely independent. This measurement count is sharp: no adaptive strategy using fewer than $d+1$ measurements guarantees identification of every single-source configuration, even when the source strength is known.
\end{theorem}

The third result concerns a fixed measurement design and establishes identification for almost every configuration of source locations. For $n\geq 1$, we define the set
\begin{equation*}
X_n = \left\{ (x_1,\ldots,x_n)\in(\R^d)^n:
x_i\neq x_k\text{ whenever }i\neq k \right\},
\end{equation*}
equipped with the restriction of $dn$-dimensional Lebesgue measure.

\begin{theorem}\label{thm:genericdesign}
Let Assumption \ref{ass:profile} hold, let $d\geq1$, and let $B\subset\R^d$ consist of $d+1$ affinely independent points. For every $n\geq1$, there exists a Lebesgue null set $E_n(B)\subset X_n$ such that, for every $(x_1,\ldots,x_n)\in X_n\setminus E_n(B)$ and every choice of positive strengths $q$, the configuration $(A,q)$ with $A=\{x_1,\ldots,x_n\}$ is uniquely identified by the traces at $B$ among all finite configurations.
The measurement count $d+1$ is sharp: no fixed set of fewer than $d+1$ measurement points guarantees identification for almost every single-source location, even when the source strength is known.
\end{theorem}

In Theorem \ref{thm:genericdesign}, the measurement design is fixed independently of the unknown configuration. For each finite cardinality, the exceptional location configurations form a Lebesgue null set that may depend on the design but not on the positive strengths. Thus any fixed set of $d+1$ affinely independent measurement points has the stated almost-everywhere guarantee for every finite cardinality. This does not conflict with Theorem \ref{thm:main}, whose universal guarantee allows no exceptional configurations.

Table \ref{tab:comparison} summarizes all the results.

\begin{table}[htbp]
\centering
\small
\setlength{\tabcolsep}{4pt}
\renewcommand{\arraystretch}{1.18}
\caption{Comparison of identification guarantees.
All guarantees hold among all finite configurations with positive strengths.}
\label{tab:comparison}
\begin{tabular}{
@{}
>{\raggedright\arraybackslash}p{\dimexpr0.27\linewidth-6pt\relax}
>{\centering\arraybackslash}p{\dimexpr0.10\linewidth-6pt\relax}
>{\raggedright\arraybackslash}p{\dimexpr0.33\linewidth-6pt\relax}
>{\raggedright\arraybackslash}p{\dimexpr0.30\linewidth-6pt\relax}
@{}
}
\toprule
\textbf{Guarantee}
&
\textbf{Sensors}
&
\textbf{Design or hypothesis}
&
\textbf{Scope}
\\
\midrule

Sharp universal threshold
(Theorem \ref{thm:main})
&
$dN+1$
&
$d\geq2$; every AGP design
&
Every configuration with at most $N$ sources
\\
\addlinespace

Sharp one-dimensional threshold
(Theorem \ref{thm:main-1d})
&
$2$
&
$d=1$; any two distinct sensors
&
Every finite configuration
\\
\addlinespace

Single source
(Corollary \ref{cor:singleton})
&
$d+1$
&
Any $d+1$ affinely independent sensors
&
Every single-source configuration
\\
\addlinespace

Clustered design
(Theorem \ref{thm:clustered})
&
$d+1$
&
Prescribed radius and squared-reference-distance gap;
sufficiently close affinely independent sensors
&
Every configuration in the prescribed class
\\
\addlinespace

Unrestricted adaptive design
(Theorem \ref{thm:adaptive})
&
$d+1$
&
First $d$ sensors affinely independent;
the $(d+1)$st chosen using the traces at those $d$ sensors
&
Every finite configuration
\\
\addlinespace

Fixed design, almost everywhere
(Theorem \ref{thm:genericdesign})
&
$d+1$
&
Any $d+1$ affinely independent sensors
&
Almost every configuration in each finite cardinality class
\\

\bottomrule
\end{tabular}
\end{table}

\section{Preliminaries on the direct problem}\label{sec:direct}

Now we describe the well-posedness of the direct problem and derive the moment identities used in the analytic--geometric equivalence. Let
\begin{equation*}
G(x,t)=(4\pi t)^{-\frac d2}e^{-\frac{|x|^2}{4t}},
\quad t>0,
\end{equation*}
be the heat kernel. The following proposition gives the Duhamel representation for a configuration $(A,q)$ \cite[Section 2.3]{Evans:2010}.

\begin{proposition}\label{prop:duhamel}
For every finite configuration $(A,q)$, problem \eqref{eq:heat} has a unique solution $u_{A,q}\in C^1([0,\infty);\calS(\R^d))$, given by
\begin{equation}\label{eq:duhamel}
u_{A,q}(x,t)
=
\int_0^t e^{s\Delta}f_{A,q}(x)\,\d s
=
\int_0^t(G(\cdot,s)*f_{A,q})(x)\,\d s.
\end{equation}
For $t>0$, $\partial_tu_{A,q}(x,t)=e^{t\Delta}f_{A,q}(x)$. Moreover, $u_{A,q}$ is smooth in $(x,t)$ for $t>0$.
\end{proposition}

The next lemma shows that the evolving profile remains radial: the contribution of a source at $a$ to a measurement at the point $b$ depends on its location only through the distance $|b-a|$. The notation $O(d)$ denotes the set of orthogonal matrices on $\mathbb{R}^d$.

\begin{lemma}\label{lem:radialconv}
For every $t>0$, the function $G(\cdot,t)*\rho$ is radial.
\end{lemma}

\begin{proof}
Let $Q\in O(d)$. Using the substitution $y=Qz$, radiality of $G(\cdot,t)$ and $\rho$, and invariance of Lebesgue measure under orthogonal transformations, we obtain
\begin{equation*}
\begin{aligned}
(G(\cdot,t)*\rho)(Qx)
&= \int_{\R^d}G(Qx-y,t)\rho(y)\,\d y = \int_{\R^d}G(Qx-Qz,t)\rho(Qz)\,\d z\\
&= \int_{\R^d}G(x-z,t)\rho(z)\,\d z = (G(\cdot,t)*\rho)(x).
\end{aligned}
\end{equation*}
Thus $G(\cdot,t)*\rho$ is invariant under $O(d)$ and depends only on $|x|$.
\end{proof}

We next derive centered radial moment identities from equality of heat traces.

\begin{lemma}\label{lem:momentidentities}
Let $(A,q)$ and $(A',q')$ be finite configurations, and set $g=f_{A',q'}-f_{A,q}$. Fix $b\in\R^d$. If $u_{A,q}(b,t)=u_{A',q'}(b,t)$ for every $t$ in a nonempty open interval $J\subset(0,\infty)$, then, for every integer $m\geq0$,
\begin{equation}\label{eq:momentzero}
\int_{\R^d}|z-b|^{2m}g(z)\,\d z=0.
\end{equation}
\end{lemma}

\begin{proof}
By Proposition \ref{prop:duhamel}, the time trace is differentiable for \(t>0\), with \(\partial_tu_{A,q}(b,t)=e^{t\Delta}f_{A,q}(b)\).
Differentiating the assumed heat trace equality on \(J\) gives $e^{t\Delta}f_{A,q}(b)=e^{t\Delta}f_{A',q'}(b)$ for $t\in J$. Thus
\[
(G(\cdot,t)*g)(b)=\int_{\R^d}G(b-z,t)g(z)\,\d z = 0, \quad t\in J,
\]
or equivalently,
$\int_{\R^d}e^{-\frac{|z-b|^2}{4t}}g(z)\,\d z=0$.
The change of variables \(s=\frac{1}{4t}\) maps \(J\) onto another
nonempty open interval \(I\subset(0,\infty)\).
Thus
\[
F(s):=\int_{\R^d}e^{-s|z-b|^2}g(z)\,\d z
\]
vanishes for every \(s\in I\). Next we show that \(F\) is holomorphic on
$\mathcal H=\{s\in\C:\Re(s)>0\}$.
Let \(K\subset\mathcal H\) be compact and put
\(\varepsilon=\min_{s\in K}\Re(s)>0\).
For every \(m\geq0\), we have
\[
\left| \frac{\partial^m}{\partial s^m} \left(e^{-s|z-b|^2}g(z)\right) \right| \leq |z-b|^{2m}e^{-\varepsilon|z-b|^2}|g(z)|.
\]
The right-hand side is integrable because \(g\in\calS(\R^d)\). Differentiation under the integral is thus valid uniformly on compact subsets of \(\mathcal H\). It follows that \(F\) is holomorphic on \(\mathcal H\), with
\[
F^{(m)}(s) = (-1)^m \int_{\R^d}|z-b|^{2m}e^{-s|z-b|^2}g(z)\,\d z.
\]
Since \(F\) vanishes on the nonempty open real interval \(I\),
the identity theorem gives \(F(s)=0\) for \(s\in\mathcal H\).
In particular, \(F^{(m)}(s)=0\) for every real \(s>0\).
Finally, for \(0<s\leq1\),
\[
|z-b|^{2m}e^{-s|z-b|^2}|g(z)|
\leq
|z-b|^{2m}|g(z)|.
\]
By the dominated convergence theorem, as \(s\downarrow0\),
\[
0 = \lim_{s\downarrow0}F^{(m)}(s) = (-1)^m \int_{\R^d}|z-b|^{2m}g(z)\,\d z.
\]
This completes the proof of the lemma.
\end{proof}

The zeroth moment shows that one heat trace uniquely determines the total source strength. When all strengths equal one, Corollary \ref{cor:totalstrength} reduces to $|A|=|A'|$.
\begin{corollary}\label{cor:totalstrength}
Under the assumptions of Lemma \ref{lem:momentidentities},
\begin{equation}\label{eq:totalstrength}
\sum_{a\in A}q_a
=
\sum_{a'\in A'}q'_{a'}.
\end{equation}
\end{corollary}
\begin{proof}
Taking $m=0$ in \eqref{eq:momentzero} and using translation invariance of Lebesgue measure, we obtain
\begin{equation*}
0
=
\int_{\R^d}(f_{A',q'}-f_{A,q})(z)\,\d z
=
M_0\left(
\sum_{a'\in A'}q'_{a'}
-
\sum_{a\in A}q_a
\right).
\end{equation*}
Since $M_0\neq0$, the assertion follows.
\end{proof}

\section{The analytic--geometric equivalence}\label{sec:equivalence}

In this section, we establish a crucial analytic--geometric equivalence between heat traces and weighted distance measures. For a finite configuration $(A,q)$ and a measurement point $b\in\R^d$, we define the weighted distance measure
\begin{equation}\label{eq:weighteddistance}
\nu_b(A,q)=\sum_{a\in A}q_a\delta_{|b-a|},
\end{equation}
where $\delta_r$ denotes the unit Dirac measure at $r\in[0,\infty)$. This is a finite positive measure satisfying
\begin{equation*}
\nu_b(A,q)(\{r\})
=
\sum_{\substack{a\in A, |b-a|=r}}q_a,
\quad r\geq0.
\end{equation*}
Thus sources equidistant from $b$ contribute their combined strength to the same support point. When all strengths equal one, $\nu_b(A,\mathbf1)$ is the counting measure associated with the distance multiset
\begin{equation*}
D_b(A):=\{\!\{|b-a|:a\in A\}\!\},
\end{equation*}
where the double braces indicate multiplicities. In particular, we have the equivalence
\begin{equation*}
\nu_b(A,\mathbf1)=\nu_b(A',\mathbf1)
\quad\Longleftrightarrow\quad
D_b(A)=D_b(A').
\end{equation*}

We first recall an elementary property of orthogonally invariant polynomials.
\begin{lemma}\label{lem:invariantpoly}
Let $p:\R^d\to\R$ be a polynomial satisfying
\begin{equation*}
p(Qv)=p(v),
\quad Q\in O(d),\quad v\in\R^d.
\end{equation*}
Then there is a one-variable polynomial $P\in\R[r]$ such that $p(v)=P(|v|^2)$.
\end{lemma}

\begin{proof}
This is the one-vector case of the first fundamental theorem for the orthogonal group \cite[Chapter II, \S9, Theorem 2.9.A]{Weyl:1946}. We give a short direct proof.
Let $e_1=(1,0,\ldots,0)$, and define $h(t)=p(te_1)$, which is a polynomial in $t$. Then reflection in the hyperplane orthogonal to $e_1$ sends $te_1$ to $-te_1$. Orthogonal invariance therefore gives $h(-t)=h(t)$, so $h$ is an even polynomial and can be written as $h(t)=P(t^2)$ for some polynomial $P$. For any $v\in\R^d$, choose $Q\in O(d)$ such that $Qv=|v|e_1$. Then we conclude
$p(v)=p(Qv)=p(|v|e_1)=P(|v|^2)$.
\end{proof}

The next lemma expresses translated radial moments as polynomials in the squared distance.

\begin{lemma}\label{lem:Pm}
Let Assumption \ref{ass:profile} hold. For every integer $m\geq0$, there is a polynomial $P_m\in\R[r]$ of degree exactly $m$ such that
\begin{equation}\label{eq:Pmidentity}
\int_{\R^d}\rho(z-a)|z-b|^{2m}\,\d z
=
P_m(|a-b|^2),
\quad a,b\in\R^d.
\end{equation}
Moreover, there exist real coefficients $c_{m,k}$ such that
\begin{equation}\label{eq:Pmtriangular}
P_m(r)=M_0r^m+\sum_{k=0}^{m-1}c_{m,k}r^k.
\end{equation}
\end{lemma}

\begin{proof}
Set $v=a-b$ and $y=z-a$. The integral in \eqref{eq:Pmidentity} becomes
\begin{equation*}
p_m(v):=\int_{\R^d}\rho(y)|y+v|^{2m}\,\d y.
\end{equation*}
Since $|y+v|^{2m}=(|v|^2+2y\cdot v+|y|^2)^m$, the factor $|y+v|^{2m}$ is a polynomial in $v$ whose coefficients are polynomials in $y$. Note that each coefficient is integrable against $|\rho(y)|\,\d y$ because $\rho$ is Schwartz. Thus $p_m$ is a polynomial in $v$ of degree at most $2m$. For $Q\in O(d)$, the substitution $y=Qw$ and radiality of $\rho$ give
\begin{align*}
p_m(Qv)
&= \int_{\R^d}\rho(y)|y+Qv|^{2m}\,\d y
=\int_{\R^d}\rho(Qw)|Qw+Qv|^{2m}\,\d w\\
&=\int_{\R^d}\rho(w)|w+v|^{2m}\,\d w=p_m(v).
\end{align*}
By Lemma \ref{lem:invariantpoly}, there is a polynomial $P_m$ such that $p_m(v)=P_m(|v|^2)$. Since $p_m$ has degree at most $2m$, $P_m$ has degree at most $m$. In the expansion of $(|v|^2+2y\cdot v+|y|^2)^m$, the only term of total degree $2m$ in $v$ is $|v|^{2m}$. Hence the homogeneous component of degree $2m$ of $p_m$ is
\begin{equation*}
\left(\int_{\R^d}\rho(y)\,\d y\right)|v|^{2m}
=
M_0|v|^{2m}.
\end{equation*}
The coefficient of $r^m$ in $P_m$ is thus  $M_0$. Since $M_0\neq0$, $P_m$ has degree exactly $m$, and the identity \eqref{eq:Pmtriangular} follows.
\end{proof}

\begin{remark}
Radiality has two roles in the analysis. First, Lemma \ref{lem:radialconv} shows that $G(\cdot,t)*\rho$ is radial, so equality of weighted distance measures implies equality of heat traces. Second, radiality ensures that the translated moment in Lemma \ref{lem:Pm} is a polynomial in $|a-b|^2$. Its triangular form permits recovery of all weighted squared-distance moments.
\end{remark}

Next we recover the weighted squared-distance moments from the radial moment identities.

\begin{lemma}\label{lem:powersums}
Let Assumption \ref{ass:profile} hold. Let $(A,q)$ and $(A',q')$ be finite configurations, and fix $b\in\R^d$. Suppose
\begin{equation*}
\int_{\R^d}|z-b|^{2m}
\bigl(f_{A',q'}(z)-f_{A,q}(z)\bigr)\,\d z=0,
\quad m\geq0.
\end{equation*}
Then
\begin{equation*}
\sum_{a\in A}q_a|a-b|^{2m}
=
\sum_{a'\in A'}q'_{a'}|a'-b|^{2m},
\quad m\geq0.
\end{equation*}
\end{lemma}

\begin{proof}
By Lemma \ref{lem:Pm}, the assumed identities give the relation
\begin{equation}\label{eq:Pmsums}
\sum_{a\in A}q_aP_m(|a-b|^2)
=
\sum_{a'\in A'}q'_{a'}P_m(|a'-b|^2),
\quad m\geq0.
\end{equation}
For every integer $m\geq0$, set
\begin{equation*}
S_m=\sum_{a\in A}q_a|a-b|^{2m}
\quad\text{and}\quad
S'_m=\sum_{a'\in A'}q'_{a'}|a'-b|^{2m},
\end{equation*}
where the zeroth powers are understood to equal one. Since $P_0=M_0\neq0$, equation \eqref{eq:Pmsums} with $m=0$ gives $S_0=S'_0$. For $m\geq1$, equations \eqref{eq:Pmtriangular} and \eqref{eq:Pmsums} imply
\begin{equation}\label{eqn:id-powersum}
M_0(S_m-S'_m)
+
\sum_{k=0}^{m-1}c_{m,k}(S_k-S'_k)=0.
\end{equation}
If $S_k=S'_k$ for every $0\leq k<m$, this identity reduces to $M_0(S_m-S'_m)=0$. Since $M_0\neq0$, induction gives $S_m=S'_m$ for every $m\geq0$.
\end{proof}

\begin{remark}
The condition $M_0\neq0$ from Assumption \ref{ass:profile} determines the total strength from the zeroth moment and gives the nonzero leading coefficient of each $P_m$. It makes the moment system \eqref{eqn:id-powersum} triangular and invertible.
\end{remark}

The weighted moments determine the distance measure through a finite Vandermonde system.

\begin{lemma}\label{lem:momentstodistances}
Let $(A,q)$ and $(A',q')$ be finite configurations, and fix $b\in\R^d$. Suppose
\begin{equation*}
\sum_{a\in A}q_a|a-b|^{2m}
=
\sum_{a'\in A'}q'_{a'}|a'-b|^{2m},
\quad m\geq0.
\end{equation*}
Then $\nu_b(A,q)=\nu_b(A',q')$.
\end{lemma}

\begin{proof}
Let $r_1,\ldots,r_L$ be the distinct squared Euclidean distances from $b$ to the points of $A\cup A'$. If $L=0$, both configurations are empty and the assertion is immediate. Otherwise, we define
\begin{equation*}
\gamma_\ell
=
\sum_{\substack{a\in A\\|a-b|^2=r_\ell}}q_a
-
\sum_{\substack{a'\in A'\\|a'-b|^2=r_\ell}}q'_{a'},
\quad \ell=1,\ldots,L.
\end{equation*}
The assumed moment identities imply
\begin{equation*}
\sum_{\ell=1}^L\gamma_\ell r_\ell^m=0,
\quad m=0,\ldots,L-1,
\end{equation*}
with $r_\ell^0=1$. The Vandermonde matrix $(r_\ell^m)_{0\leq m\leq L-1,\,1\leq\ell\leq L}$ is invertible since the $r_\ell$ are distinct. Hence $\gamma_\ell=0$ for every $\ell$. Since $\gamma_\ell$ is the difference between the masses of $\nu_b(A,q)$ and $\nu_b(A',q')$ at $\sqrt{r_\ell}$, and both measures are supported on these distances, the two measures are equal.
\end{proof}

Conversely, equality of weighted distance measures implies equality of heat traces by radiality of the evolving profile.

\begin{lemma}\label{lem:distancestotraces}
Let $(A,q)$ and $(A',q')$ be finite configurations, and fix $b\in\R^d$. If $\nu_b(A,q)=\nu_b(A',q')$, then $u_{A,q}(b,t)=u_{A',q'}(b,t)$ for every $t>0$.
\end{lemma}

\begin{proof}
By Lemma \ref{lem:radialconv}, for each $s>0$, there is a function $\psi_s:[0,\infty)\to\R$ such that $(G(\cdot,s)*\rho)(x)=\psi_s(|x|)$. Consequently,
\begin{equation*}
\begin{aligned}
e^{s\Delta}f_{A,q}(b)
=\sum_{a\in A}q_a(G(\cdot,s)*\rho)(b-a)=\sum_{a\in A}q_a\psi_s(|b-a|)=\int_{[0,\infty)}\psi_s(r)\,\nu_b(A,q)(\d r).
\end{aligned}
\end{equation*}
The formula holds also for $(A',q')$. Equality of the distance measures then gives $e^{s\Delta}f_{A,q}(b)=e^{s\Delta}f_{A',q'}(b)$ for every $s>0$. Integrating over $s\in(0,t)$ and using Proposition \ref{prop:duhamel} yields the desired assertion.
\end{proof}

The preceding results give the following analytic--geometric equivalence.

\begin{theorem}\label{thm:equivalence}
Let Assumption \ref{ass:profile} hold. Let $(A,q)$ and $(A',q')$ be finite configurations, let $b\in\R^d$, and let $J\subset(0,\infty)$ be a nonempty open interval. Then the following statements are equivalent:
\begin{enumerate}
\item[{\rm(i)}]
$u_{A,q}(b,t)=u_{A',q'}(b,t)$ for every $t\in J$;
\item[{\rm(ii)}]
$\nu_b(A,q)=\nu_b(A',q')$ as measures on $[0,\infty)$;
\item[{\rm(iii)}]
$u_{A,q}(b,t)=u_{A',q'}(b,t)$ for every $t>0$.
\end{enumerate}
\end{theorem}

\begin{proof}
Assume {\rm(i)}. Lemma \ref{lem:momentidentities} gives all centered radial moment identities for $g=f_{A',q'}-f_{A,q}$. Lemma \ref{lem:powersums} then gives equality of all weighted squared-distance moments, and Lemma \ref{lem:momentstodistances} yields {\rm(ii)}. Lemma \ref{lem:distancestotraces} proves that {\rm(ii)} implies {\rm(iii)}. Finally, {\rm(iii)} immediately implies {\rm(i)}.
\end{proof}

\begin{remark}
Theorem \ref{thm:equivalence} requires neither AGP nor a bound on the number of source locations. It shows that one heat trace determines the distances carrying nonzero mass and the total strength at each such distance, but does not distinguish sources at the same distance.
\end{remark}

\section{Universal bounds}\label{sec:bounds}

In this section, we prove Theorems \ref{thm:main} and \ref{thm:main-1d}, and Corollary \ref{cor:singleton}. We use the following notation frequently. For finite $(A,q)$ and $(A',q')$, we define the witness set
\begin{equation*}
\calW((A,q),(A',q')) = \{b\in\R^d:\nu_b(A,q)=\nu_b(A',q')\}.
\end{equation*}
By Theorem \ref{thm:equivalence}, a point belongs to the set $\calW((A,q),(A',q'))$ precisely when the two configurations produce the same heat trace there on one, and hence every, nonempty open positive-time interval. For unit-strength configurations, we simplify the notation
\begin{equation*}
\calW(A,A')
:=
\calW((A,\mathbf1),(A',\mathbf1))
=
\{b\in\R^d:D_b(A)=D_b(A')\}.
\end{equation*}

For $a,a'\in\R^d$, $a\neq a'$, we define the perpendicular
bisector hyperplane
\begin{equation*}
H(a,a')
=
\calW(\{a\},\{a'\})
=
\{x\in\R^d:|x-a|=|x-a'|\}.
\end{equation*}
Squaring the equality and cancelling terms gives
\begin{equation*}
2x\cdot(a'-a)=|a'|^2-|a|^2.
\end{equation*}
Since $a'-a\neq0$, this is a nontrivial affine linear equation.
Hence $H(a,a')$ is a proper affine hyperplane.

\subsection{Tools for deriving upper and lower bounds}

Now we develop several tools for deriving upper and lower bounds.
We first establish the geometric upper bound.
\begin{lemma}\label{lem:geometricupper}
Let $d\geq1$ and $N\geq1$. Let $(A,q)$ and $(A',q')$ be
finite configurations with positive strengths and $|A|\leq N$.
Let $B\subset\R^d$ be in AGP and suppose
\begin{equation*}
B\subset\calW((A,q),(A',q')).
\end{equation*}
If the two configurations are distinct, then $|B|\leq dN$.
Equivalently, if $|B|\geq dN+1$, then $A=A'$ and
$q_a=q'_a$ for every $a\in A$.
No cardinality bound is imposed on $A'$.
\end{lemma}

\begin{proof}
Suppose $|B|\geq dN+1$. By the definition of the witness set $\calW((A,q),(A',q'))$, we have
\begin{equation}\label{eq:weighteddistanceagreement}
\nu_b(A,q)=\nu_b(A',q'),
\quad b\in B.
\end{equation}
If $A=\varnothing$, then equality at any point of the nonempty set $B$ and positivity of strengths imply $A'=\varnothing$. Hence, it suffices to consider the case that $1\leq |A|=:n\leq N$.
First we show that $A'\subseteq A$. Suppose instead that $a'\in A'\setminus A$. For every $b\in B$, the distance $|b-a'|$ carries positive mass in $\nu_b(A',q')$. By \eqref{eq:weighteddistanceagreement}, there is
some $a\in A$ such that $|b-a'|=|b-a|$. Consequently,
\begin{equation}\label{eq:weightedbisectorcover}
B\subseteq\bigcup_{a\in A}H(a',a).
\end{equation}
Since $a'\notin A$, each $H(a',a)$ is a proper affine hyperplane. Then,
AGP implies $|B\cap H(a',a)|\leq d$, and hence $|B|\leq d|A|\leq dN$,
which contradicts the assumption $|B|\ge dN+1$. Thus $A'\subseteq A$, and in particular $|A'|\leq N$.
Interchanging the two configurations and repeating the
preceding argument gives $A\subseteq A'$. Hence $A=A'$.
It remains to identify the strengths. Fix $a\in A$. By AGP, the union $\bigcup_{c\in A\setminus\{a\}}H(a,c)$
contains at most $d(|A|-1)$ points of $B$. Since $|B|\geq dN+1>d(|A|-1)$, there exists $b\in B$
such that $|b-a|\neq |b-c|$ for every $c\in A\setminus\{a\}$.
Then equality of the distance measures at the point $b$ yields
\begin{equation*}
q_a = \nu_b(A,q)(\{|b-a|\}) = \nu_b(A',q')(\{|b-a|\}) = q'_a.
\end{equation*}
Since $a$ can be chosen arbitrarily, all strengths agree. This completes the proof of the lemma.
\end{proof}

\begin{remark}
The analytic--geometric equivalence at a single point does
not require AGP. AGP enters only through the fact that a
proper affine hyperplane contains at most $d$ measurement
points. More generally, if
$\sup_H|B\cap H|\leq K$,
where the supremum is over proper affine hyperplanes and
$K\in\mathbb{N}$, the same proof shows that
$KN+1$ measurements suffice. AGP is the case $K=d$.
\end{remark}

\begin{remark}
Positivity of strengths is essential to the support argument above:
each source contributes positive mass at its distance
from every measurement point. With strengths of mixed
signs, contributions from equidistant sources may cancel,
so the argument no longer applies.
\end{remark}

Combining the geometric bound in Lemma \ref{lem:geometricupper} with the analytic-geometric equivalence in Theorem
\ref{thm:equivalence} gives the upper bound
for heat-trace measurements.

\begin{corollary}\label{cor:analyticupper}
Let Assumption \ref{ass:profile} hold, let $d\geq1$ and
$N\geq1$, and let $B\subset\R^d$ be in AGP with
$|B|\geq dN+1$. Let $(A,q)$ and $(A',q')$ be finite
configurations with $|A|\leq N$. If
\begin{equation*}
u_{A,q}(b,t)=u_{A',q'}(b,t)
\quad\text{for all}\quad b\in B~~\text{and}~~ t\in J,
\end{equation*}
then $A=A'$ and $q_a=q'_a$ for every $a\in A$.
\end{corollary}

\begin{proof}
Theorem \ref{thm:equivalence} gives
$B\subset\calW((A,q),(A',q'))$.
Lemma \ref{lem:geometricupper} then gives the assertion.
\end{proof}

For sharpness, it suffices to construct ambiguities
with unit strengths.
The lower bound requires selecting $d$ points on each of
$N$ witness hyperplanes while preserving AGP. First, the
next result provides a generic selection of points on
distinct hyperplanes.

\begin{lemma}\label{lem:genericselection}
Let $d\geq2$, and let $H_1,\ldots,H_N$ be distinct affine
hyperplanes in $\R^d$. Then there is a finite AGP set
$B\subset\R^d$ such that
\begin{equation*}
B\subset\bigcup_{j=1}^NH_j
\quad\mbox{and}\quad
|B\cap H_j|=d
\quad(j=1,\ldots,N).
\end{equation*}
In particular, $|B|=dN$.
\end{lemma}

\begin{proof}
We select the points successively, assigning $d$ points
to each $H_j$, while maintaining AGP and ensuring that
each selected point lies on no other $H_k$. Let $S$ be
the set of points already selected, and suppose that
fewer than $d$ points have been assigned to $H_j$. We choose the next point in $H_j$ outside
\begin{equation*}
\bigcup_{k\ne j}(H_j\cap H_k)
\;\cup\!
\bigcup_{\substack{T\subseteq S\\1\le |T|\le d}}
\bigl(H_j\cap\aff(T)\bigr).
\end{equation*}
Every nonempty intersection above is a proper affine
subspace of $H_j$. Indeed, the $H_k$ are distinct.
Moreover, $H_j\subseteq\aff(T)$ (with $\aff(T)$ being the affine hull) would force $|T|=d$
and $\aff(T)=H_j$, contradicting the fact that $S$
contains fewer than $d$ points on $H_j$.
Since $\dim H_j=d-1\geq1$, these finitely many proper
affine subspaces have zero $(d-1)$-dimensional Lebesgue
measure in $H_j$, and hence cannot cover it.
The new point therefore exists and preserves both
properties. Repeating this procedure until each $H_j$
has $d$ assigned points yields an AGP set
$B\subseteq\bigcup_{j=1}^N H_j$ with $|B\cap H_j|=d$
for every $j$ and $|B|=dN$.
\end{proof}

\begin{remark}
Excluding $H_j\cap H_k$ ensures that each selected point
lies only on its assigned hyperplane. Without this
restriction, points assigned to other hyperplanes could
also lie on $H_j$, so assigning $d$ points to each
hyperplane could place more than $d$ selected points
on $H_j$ and violate AGP. This is relevant to the
construction below, where the witness hyperplanes share
a $(d-2)$-dimensional affine subspace.
\end{remark}

Next we construct two interlaced regular polygons, whose unit-strength configurations are indistinguishable on $N$ hyperplanes. It is used in establishing the lower bound.

\begin{lemma}\label{lem:reflection}
Let $d\geq2$ and $N\geq1$.
There exist distinct $N$-point sets $A,A'\subset\R^d$
and $N$ distinct affine hyperplanes $H_0,\ldots,H_{N-1}$
such that
\begin{equation*}
H_m\subset\calW(A,A'),
\quad m=0,\ldots,N-1.
\end{equation*}
\end{lemma}

\begin{proof}
We embed $\R^2$ into $\R^d$ as
$\Pi=\R^2\times\{0\}^{d-2}$, and identify $\R^2$ with $\mathbb{C}$. Let
\begin{align*}
A
=
\left\{e^{\frac{2ki\pi }{N}}:k=0,\ldots,N-1\right\}\subset\Pi~~\text{and}~~
A' = e^{\frac{i\pi}{N}}A = \left\{e^{\frac{i(2k+1)\pi}{N}}:k=0,\ldots,N-1\right\}\subset\Pi.
\end{align*}
Both sets have exactly \(N\) points.
They are distinct.
Indeed, if \(A'=A\), rotation by \(\frac{\pi}{N}\) would preserve \(A\),
and for some integer \(k\), $\frac{\pi}{N}\equiv\frac{2\pi k}{N}\pmod{2\pi}$.
Multiplication by \(\frac{N}{\pi}\) gives
\(1\equiv2k\pmod{2N}\), which is impossible. For \(N=1\), the sets are \(A=\{1\}\) and \(A'=\{-1\}\).
Now for $m=0,\ldots,N-1$, we define
\begin{equation*}
\theta_m=\frac{\pi}{2N}+\frac{m\pi}{N},
\quad
\ell_m=\{re^{i\theta_m}:r\in\R\}.
\end{equation*}
These unoriented lines are distinct. Indeed,
$\ell_m=\ell_{m'}$ would imply
\begin{equation*}
\theta_m-\theta_{m'}
=
\frac{(m-m')\pi}{N}
\in\pi\mathbb Z,
\end{equation*}
so $m-m'$ is a multiple of $N$. The stated index range then gives $m=m'$. Then reflection in $\ell_m$ sends polar angle $\alpha$ to $2\theta_m-\alpha$. For the point of $A$ with angle
$\frac{2\pi k}{N}$, the reflected angle is
\begin{equation*}
2\theta_m-\frac{2\pi k}{N} = \frac{(2m-2k+1)\pi}{N}.
\end{equation*}
This is an odd multiple of $\frac{\pi}{N}$, so the reflected
point belongs to $A'$. As $k$ varies modulo $N$, so
does $m-k$, and reflection maps $A$ bijectively onto
$A'$. Now we set $H_m=\ell_m\times\R^{d-2}$.
These are distinct affine hyperplanes in $\R^d$.
Let $R_m$ denote orthogonal reflection across $H_m$.
The planar calculation gives $R_m(A)=A'$.
For any $b\in H_m$, we have $R_mb=b$, and
\begin{equation*}
|b-a|
=
|R_mb-R_ma|
=
|b-R_ma|.
\end{equation*}
Meanwhile, as $a$ ranges over $A$, $R_ma$ ranges bijectively over
$A'$. Therefore $D_b(A)=D_b(A')$. Since $b\in H_m$
can be chosen arbitrarily, $H_m\subset\calW(A,A')$.
\end{proof}

The next result gives the lower-bound construction.

\begin{theorem}\label{thm:lower}
Let Assumption \ref{ass:profile} hold, and let
$d\geq2$ and $N\geq1$.
There exist distinct $N$-point sets $A,A'\subset\R^d$
and an AGP set $B\subset\R^d$, $|B|=dN$, such that
\begin{equation*}
u_{A,\mathbf1}(b,t)=u_{A',\mathbf1}(b,t),
\quad b\in B,\quad t>0.
\end{equation*}
\end{theorem}

\begin{proof}
Lemma \ref{lem:reflection} supplies distinct $N$-point
sets $A,A'$ and distinct hyperplanes $H_0,\ldots,H_{N-1}$
such that $H_m\subset\calW(A,A')$.
Lemma \ref{lem:genericselection} gives an AGP set
\begin{equation*}
B\subset\bigcup_{m=0}^{N-1}H_m,
\end{equation*}
with exactly $d$ points on each hyperplane and no
selected point on two different witness hyperplanes.
Thus $|B|=dN$ and $B\subset\calW(A,A')$.
Since equality of distance multisets is equivalent
to equality of the unit-strength distance measures,
Theorem \ref{thm:equivalence} gives the asserted equality
of heat traces for every $t>0$.
\end{proof}


\begin{remark}
The upper and lower bounds use the hyperplane
occupancy count in different ways. In the upper bound,
an additional competing location would force the
measurement set into a union of at most $N$ proper
bisector hyperplanes. AGP permits at most $d$ measurement
points on each, giving the bound $dN$. Once the supports
coincide, the same geometry supplies measurements that
separate each source from all the others and determine
its strength. In the lower bound, unit-strength
configurations have $N$ witness hyperplanes, and one
can select $d$ points on each while preserving AGP.
Thus the sharpness obstruction persists even when all
strengths are known.
\end{remark}

\subsection{Proof of Theorems \ref{thm:main} and
\ref{thm:main-1d}}\label{sec:proofs}

Now we prove the main results.

\begin{proof}[Proof of Theorem \ref{thm:main}]
Part {\rm(i)} follows directly from Corollary
\ref{cor:analyticupper}.
Theorem \ref{thm:lower} proves part {\rm(ii)} using
unit-strength configurations. Every subset of the
constructed measurement set remains in AGP and gives
the same ambiguity. Therefore no positive measurement
count at most $dN$ guarantees identification for every
AGP design. Together with part {\rm(i)}, this proves
that $dN+1$ is the sharp universal threshold.
\end{proof}

The following example illustrates the sharpness of the
universal threshold for $d=2$ and $N=2$: four measurement
points in AGP may fail to distinguish two source
configurations. It also illustrates the labeling
ambiguity across measurement points.

\begin{proposition}\label{prop:4pointcounterexample}
Let $d=2$, and suppose that $\rho$ satisfies
Assumption \ref{ass:profile}. Set
\begin{equation*}
A=\{(1,0),(-1,0)\},
\quad
A'=\{(0,1),(0,-1)\},
\quad
B=\{(1,1),(2,2),(1,-1),(2,-2)\}.
\end{equation*}
Then $A\neq A'$, $B$ is in AGP, and
\begin{equation*}
u_{A,\mathbf1}(b,t)=u_{A',\mathbf1}(b,t),
\quad b\in B,\quad t>0.
\end{equation*}
Thus four measurement points in AGP do not in general
uniquely identify two source locations
even when all strengths are known and equal to one.
\end{proposition}

\begin{proof}
The set $B$ consists of two points on each of the
parallel lines $x_1=1$ and $x_1=2$. Any three points
of $B$ contain two on one line and one on the other,
so they are noncollinear. Hence $B$ is in AGP.
Besides, direct calculation gives
\begin{equation*}
\begin{array}{c|c|c}
b&
\{\!\{|b-a|^2:a\in A\}\!\}&
\{\!\{|b-a'|^2:a'\in A'\}\!\}\\ \hline
(1,1)&\{\!\{1,5\}\!\}&\{\!\{1,5\}\!\}\\
(2,2)&\{\!\{5,13\}\!\}&\{\!\{5,13\}\!\}\\
(1,-1)&\{\!\{1,5\}\!\}&\{\!\{1,5\}\!\}\\
(2,-2)&\{\!\{5,13\}\!\}&\{\!\{5,13\}\!\}
\end{array}
\end{equation*}
Thus $D_b(A)=D_b(A')$ for every $b\in B$, and Theorem \ref{thm:equivalence} yields equality of the unit-strength heat traces for all $t>0$,
although $A\neq A'$.
\end{proof}

\begin{remark}
This example illustrates the distinction between labeled distances and distance multisets. For unit-strength configurations, the heat trace at each measurement point determines the distances as an unordered multiset, but does not specify how to match them across measurement points. Reflection across $x_1=x_2$ maps $A$ onto $A'$ and fixes $(1,1)$ and $(2,2)$, whereas reflection across $x_1=-x_2$ also maps $A$ onto $A'$ and fixes $(1,-1)$ and $(2,-2)$. Thus, the two configurations give the same distance multiset at every measurement point, even though their source locations are different.
\end{remark}

The next result shows that an \textsl{a priori} restriction on the source region also permits uniqueness from measurements confined to a hyperplane.

\begin{corollary}\label{cor:hyperplane-region}
Let Assumption \ref{ass:profile} hold, and let $d\geq2$ and $N\geq1$. Let $D\subseteq\R^d$, let $P$ be an
affine hyperplane, and denote reflection across $P$ by $R_P$. Suppose that $ D\cap R_P(D)\subseteq P$. Let $B\subseteq P$ consist of $(d-1)N+1$ points in AGP within $P$ such that every subset of at most $d$ points is affinely independent. Let $(A,q)$ and $(A',q')$ be finite configurations with $A,A'\subseteq D$ and $|A|\leq N$.
If
\begin{equation*}
u_{A,q}(b,t)=u_{A',q'}(b,t), \quad b\in B,\quad t\in J,
\end{equation*}
then $A=A'$ and $q_a=q'_a$ for every $a\in A$. In particular, for any bounded set $D$, such a design can be obtained by choosing $P$ so that $D$ lies entirely in one of its open half-spaces.
\end{corollary}

\begin{proof}
Theorem \ref{thm:equivalence} gives equality of the weighted distance measures at every point of $B$.
For distinct $x,y\in D$, the equality $H(x,y)=P$ would imply $y=R_P(x)$. Then $x\in D\cap R_P(D)\subseteq P$, which would give $y=x$, a contradiction. Thus every nonempty
intersection $P\cap H(x,y)$ is a proper affine subspace of $P$ and contains at most $d-1$ points of $B$.
Suppose $x\in A'\setminus A$. Positivity and equality of the distance measures imply
\begin{equation*}
B\subseteq\bigcup_{a\in A}\bigl(P\cap H(x,a)\bigr).
\end{equation*}
Consequently, we have
\begin{equation*}
(d-1)N+1=|B|\leq(d-1)|A| \leq(d-1)N.
\end{equation*}
This leads to a contradiction. Hence $A'\subseteq A$ and $|A'|\leq N$. Interchanging the configurations gives $A\subseteq A'$, so $A=A'$. For each $a\in A$, the set $B\cap\bigcup_{c\in A\setminus\{a\}}H(a,c)$ has at most $(d-1)(|A|-1)$ points. Thus there is a measurement point at which $a$ has a different
distance from every other source. Equality of the
mass at this distance gives $q_a=q'_a$.
Finally, if $D$ lies in an open half-space bounded
by $P$, then $D\cap R_P(D)=\varnothing$. Then the result follows immediately.
\end{proof}
\begin{proof}[Proof of Theorem \ref{thm:main-1d}] By Theorem \ref{thm:equivalence}, we have
$\nu_{b_j}(A,q)=\nu_{b_j}(A',q')$ for  $j=1,2$. Next we extend $q$ and $q'$ to functions on $\R$ by setting
them equal to zero outside $A$ and $A'$, respectively, and define the finitely supported function
\begin{equation*}
h(x)=q(x)-q'(x), \quad x\in\R.
\end{equation*}
Fix $j\in\{1,2\}$. Then equality of the masses at every
distance $r>0$ gives
\begin{equation*}
h(b_j+r)+h(b_j-r)=0,
\end{equation*}
while equality at distance zero gives $h(b_j)=0$. Thus
\begin{equation*}
h(2b_j-x)=-h(x), \quad x\in\R.
\end{equation*}
Composing the two reflections yields
\begin{equation*}
h\bigl(x+2(b_2-b_1)\bigr) = h\bigl(2b_2-(2b_1-x)\bigr) = -h(2b_1-x) = h(x).
\end{equation*}
Since $b_1\neq b_2$, the number $L=2(b_2-b_1)$ is nonzero. If $h(x_0)\neq0$, then
\begin{equation*}
h(x_0+kL)=h(x_0)\neq0, \quad k=0,1,2,\ldots,
\end{equation*}
which contradicts the finite support condition. Hence $h=0$. Then positivity of the original strengths gives $A=A'$ and $q_a=q'_a$ for every $a\in A$. For sharpness, fix one measurement point $b\in\R$
and choose $r>0$. The distinct singleton configurations at $b+r$ and $b-r$, both with unit strength, have the same distance measure at $b$.
By Theorem \ref{thm:equivalence}, their heat traces there agree for every $t>0$. Thus one measurement cannot guarantee identification, and the sharp
count is two.
\end{proof}
\begin{proof}[Proof of Corollary \ref{cor:singleton}]
For $d\geq2$, affine independence implies that $\{b_0,\ldots,b_d\}$ is in AGP, so identification among all finite configurations follows from Theorem \ref{thm:main}{\rm(i)} with $N=1$.
For $d=1$, the same conclusion follows from Theorem \ref{thm:main-1d}. For necessity, any nonempty set of fewer than $d+1$ measurement points is contained in an affine
hyperplane $H$. Now we choose $a\notin H$, and let $a'$ be its reflection across $H$. Then $a\neq a'$, but their distances to every measurement point agree. Theorem \ref{thm:equivalence} then gives identical traces for the two singleton configurations with unit strengths.
\end{proof}

\section{Design-dependent and adaptive measurement bounds}
\label{sec:design}

In this section, we present the proofs of the design-dependent and adaptive measurement bounds in Theorems \ref{thm:clustered}, \ref{thm:adaptive} and \ref{thm:genericdesign}. In one dimension, the sharper conclusion of Theorem \ref{thm:main-1d} already applies.
First we prove Theorem \ref{thm:clustered} on clustered measurements with separated reference distances.

\begin{proof}[Proof of Theorem \ref{thm:clustered}]
By Theorem \ref{thm:equivalence}, the weighted distance measures agree at every measurement point.
If $A=\varnothing$, positivity gives $A'=\varnothing$, and the conclusion follows. Hence assume that $n:=|A|\geq1$.
By the definition of the set $\mathfrak A(R,\delta;b_0)$, the squared reference distances are distinct.
Thus we can label the points of \(A\) so that
\[
r_1<r_2<\cdots<r_n, \quad\mbox{with } r_i=|a_i-b_0|^2.
\]
(Clearly, for $n=1$, we have the unique label.) Equality of the measures at $b_0$ partitions $A'$ into nonempty
groups
\[
A'_i=\{a'\in A':|a'-b_0|^2=r_i\}, \quad \sum_{a'\in A'_i}q'_{a'}=q_{a_i}, \quad i=1,\ldots,n.
\]
In particular, both configurations lie in $\overline{B(b_0,R)}$. Fix $j\in\{1,\ldots,d\}$ and set $h_j=b_j-b_0$.
For every $x\in\{a_i\}\cup A'_i$,
\[
|x-b_j|^2=r_i+|h_j|^2-2h_j\cdot(x-b_0).
\]
Hence these squared distances belong to the interval
\[
I_{i,j} = \bigl[r_i+|h_j|^2-2R|h_j|,\, r_i+|h_j|^2+2R|h_j|\bigr].
\]
For each fixed $j$, these intervals are pairwise disjoint, since
$$r_{i+1}-r_i\geq\delta>4R|h_j|,\quad\mbox{for }i=1,\ldots,n-1.$$
For $n=1$, this condition is vacuous. Let $a'\in A'_i$. Now positivity and equality of the weighted distance measures at $b_j$ imply, for some $k$,
$|a'-b_j|^2=|a_k-b_j|^2$. The common squared distance belongs to both $I_{i,j}$ and $I_{k,j}$, so disjointness forces $k=i$. Since $j$ was arbitrary, and the reference distances also agree,
\[
|a'-b_j|^2=|a_i-b_j|^2,
\quad j=0,\ldots,d.
\]
Subtracting the equation at $b_0$ from the others gives
$$(b_j-b_0)\cdot(a'-a_i)=0, \quad  j=1,\ldots,d.$$
Affine independence implies $a'=a_i$.
Thus $A'_i=\{a_i\}$, and equality of the total strengths in
this group gives $q'_{a_i}=q_{a_i}$.
\end{proof}

\begin{remark}\label{rem:clustereddesign}
Let $e_1,\ldots,e_d$ be the standard basis vectors.
The points $b_j=b_0+\varepsilon e_j$, $j=1,\ldots,d$
satisfy the measurement conditions of Theorem \ref{thm:clustered}
whenever $0<\varepsilon<\frac{\delta}{4R}$.
\end{remark}

\begin{remark}\label{rem:labelledreconstruction}
Once the distances have been consistently labeled, the locations
can be recovered by solving a linear system.
Set $h_j=b_j-b_0$, and let $H$ be the $d\times d$ matrix whose
$j$-th row is $h_j^{\mathsf T}$.
With $r_{i,j}=|a_i-b_j|^2$, we have
\[
2h_j\cdot(a_i-b_0)
=
|h_j|^2+r_{i,0}-r_{i,j},
\quad j=1,\ldots,d,
\]
and hence
\[
a_i=b_0+\frac12H^{-1}
\begin{pmatrix}
|h_1|^2+r_{i,0}-r_{i,1}\\
\vdots\\
|h_d|^2+r_{i,0}-r_{i,d}
\end{pmatrix}.
\]
The strength is recovered from the weighted distance measure at $b_0$:
$q_{a_i}=\nu_{b_0}(A,q)(\{\sqrt{r_{i,0}}\})$.
\end{remark}

The separation parameters can be read from the first measurement
whenever the reference distances are distinct, leading to an
adaptive identification from a separating reference point.

\begin{corollary}\label{cor:adaptiveclustered}
Let Assumption \ref{ass:profile} hold.
Let $(A,q)$ be a finite configuration, and let $b_0\in\R^d$ satisfy
\[
|a-b_0|\neq|\widetilde a-b_0|
\quad\text{for all distinct }a,\widetilde a\in A.
\]
Then $(A,q)$ can be uniquely identified among all finite
configurations using $d+1$ measurements, with the remaining
points chosen after observing the trace at $b_0$.
For every fixed finite set $A$, almost every reference point $b_0$
satisfies this separation condition.
\end{corollary}
\begin{proof}
By Theorem \ref{thm:equivalence}, the heat trace at $b_0$ uniquely
determines $\nu_{b_0}(A,q)$.
Under the separation hypothesis, its support has cardinality $n=|A|$.
If $n=0$, positivity identifies the empty configuration.
If $n=1$, choose any $d$ further points forming an affine simplex
with $b_0$.
Corollary \ref{cor:singleton} gives uniqueness among all finite
configurations.
Suppose $n\geq2$.
The support of the reference weighted distance measure determines
\[
R=\max_{a\in A}|a-b_0|,
\quad
\delta=
\min_{\substack{a,\widetilde a\in A\\a\neq\widetilde a}}
\bigl||a-b_0|^2-|\widetilde a-b_0|^2\bigr|.
\]
The hypothesis gives $\delta>0$, and $n\geq2$ gives $R>0$. Now we choose $b_j=b_0+\varepsilon e_j$, $j=1,\ldots,d$, with $0<\varepsilon<\frac{\delta}{4R}$. Theorem \ref{thm:clustered} gives uniqueness among all finite configurations.
Finally, for each distinct pair $a,\widetilde a\in A$, the reference points where their distances coincide form the proper bisector hyperplane $H(a,\widetilde a)$.
Thus the exceptional reference points form a finite union of proper hyperplanes, which has Lebesgue measure zero.
\end{proof}

Next, we prove Theorem \ref{thm:adaptive}.
\begin{proof}[Proof of Theorem \ref{thm:adaptive}]
In the proof, we first show that the traces at the first $d$ affinely independent measurement points restrict all competing source locations to a finite candidate set. Then we choose one additional measurement point so that its distances to these candidate locations
are pairwise distinct, which allows the source locations and their strengths to be uniquely recovered. The proof is divided into three steps.\\
\textbf{Step 1. Restrict the competing source to a finite candidate set.} Fix affinely independent $b_0,\ldots,b_{d-1}$.
By Theorem \ref{thm:equivalence}, the heat traces at $b_0,\ldots,b_{d-1}$ determine $\nu_{b_j}(A,q)$ for $j=0,\ldots,d-1$.
If $A=\varnothing$, positivity of the strengths already gives uniqueness.
Hence assume that $A\neq\varnothing$. Then positivity of the strengths shows that these measures determine the sets
\[
\mathcal R_j=\{|a-b_j|^2:a\in A\},
\quad j=0,\ldots,d-1.
\]
Now we define the candidate location set
\[
S=\{x\in\R^d:|x-b_j|^2\in\mathcal R_j
\text{ for every }j=0,\ldots,d-1\}.
\]
By Theorem \ref{thm:equivalence} and positivity of the strengths, every finite configuration $(A',q')$ reproducing the traces at $b_0,\ldots,b_{d-1}$ satisfies $A'\subseteq S$. In particular, $A\subseteq S$.
We first show that $S$ is finite. For each tuple $(r_0,\ldots,r_{d-1})\in\mathcal R_0\times\cdots\times\mathcal R_{d-1}$, consider $|x-b_j|^2=r_j$ for all $j=0,\ldots,d-1$. Subtracting the equation at $b_0$ from the others gives
\begin{equation*}
2(b_j-b_0)\cdot x = |b_j|^2-|b_0|^2+r_0-r_j,
\quad j=1,\ldots,d-1.
\end{equation*}
Since $b_0,\ldots,b_{d-1}$ are affinely independent, the $d-1$ linear equations are independent, so their solution set is a one-dimensional affine subspace of $\R^d$. Intersecting this line with $|x-b_0|^2=r_0$ gives at most
two points. For $d=1$, there are no linear equations, and $|x-b_0|^2=r_0$  has at most two solutions. Thus $S$ is determined by the traces at $b_0,\ldots,b_{d-1}$ and satisfies $|S|\leq 2\prod_{j=0}^{d-1}|\mathcal R_j|$.\\
\textbf{Step 2. Identify locations and strengths using an additional point.} To distinguish all candidate locations, choose the $(d+1)$st measurement point so that
\[
b_d\notin
\aff\{b_0,\ldots,b_{d-1}\}
\;\cup\;
\bigcup_{\substack{x,y\in S\\x\neq y}}H(x,y).
\]
Such a choice is possible because the set on the right is a finite union of proper affine hyperplanes and hence cannot cover $\R^d$. Moreover, it depends only on $b_0,\ldots,b_{d-1}$ and their traces. This choice makes $b_0,\ldots,b_d$ affinely independent and the distances from $b_d$ to the points of $S$ pairwise distinct.
Now let $(A',q')$ be any finite configuration with positive strengths that reproduces the traces at all points $b_0,\ldots,b_d$. The observations at $b_0,\ldots,b_{d-1}$ imply $A'\subseteq S$. We extend $q$ and $q'$ to $S$ by setting them equal to zero outside $A$ and $A'$, respectively. By Theorem \ref{thm:equivalence}, equality of the traces at $b_d$ gives
$\nu_{b_d}(A,q)=\nu_{b_d}(A',q')$.
Since the distances from $b_d$ to the points of $S$ are pairwise distinct, we obtain
\[
q_x
=
\nu_{b_d}(A,q)(\{|b_d-x|\})
=
\nu_{b_d}(A',q')(\{|b_d-x|\})
=
q'_x,
\quad x\in S.
\]
Strict positivity of the strengths on their respective
source sets now yields $A=A'$ and $q=q'$.\\
\textbf{Step 3. Prove sharpness}.
Fix any adaptive strategy using at most $d$ measurements, and restrict our attention to single-source configurations with unit strength. By radiality and Theorem~\ref{thm:equivalence}, the full
trace at each measurement point is uniquely determined by, and
uniquely determines, the distance to the source. Let $b_0$ be the first measurement point. All source locations on the unit sphere centered at $b_0$ produce the same trace at $b_0$.  We proceed inductively. Let $C_0=\{b_0+v:|v|=1\}$. Let $b_j$ be the next measurement point selected by the strategy based on the common traces at $b_0,\ldots,b_{j-1}$, and define
\begin{equation*}
C_j=
\left\{
w=b_0+v:\ |v|=1,\quad
v\cdot(b_i-b_0)=0,\quad i=1,\ldots,j
\right\}.
\end{equation*}
Since $j\leq d-1$, $\operatorname{span}\{b_i-b_0:1\leq i\leq j\}^\perp$ has dimension at least one, so $C_j$ contains at least two points.
Moreover, $C_j\subseteq C_{j-1}$, so every source location in $C_j$
reproduces all preceding traces. For every $w = b_0+v \in C_j$,
\begin{equation*}
|w-b_j|^2 = |b_0+v-b_j|^2=1+|b_j-b_0|^2.
\end{equation*}
Therefore, all unit-strength single sources with locations in
$C_j$ produce the same trace at $b_k$ for all $0\le k\le j$. When the strategy stops after at most $d$ measurements, at least two distinct single-source configurations remain indistinguishable. This proves the claimed sharpness.
\end{proof}

Last, we prove the fixed-design result in Theorem \ref{thm:genericdesign} by constructing a Lebesgue null exceptional set using polynomial compatibility conditions for distances to $d+1$ affinely independent measurement points.

\begin{proof}[Proof of Theorem \ref{thm:genericdesign}]
Let $B=\{b_0,\ldots,b_d\}$. Fix $n\geq1$ and view  $x=(x_1,\ldots,x_n)$ as a variable. For each tuple $\alpha=(i_0,\ldots,i_d)\in\{1,\ldots,n\}^{d+1}$, we define
\begin{equation}\label{eqn:sq-dist}
r_j^\alpha(x)=|x_{i_j}-b_j|^2, \quad j=0,\ldots,d.
\end{equation}
The proof is lengthy and is divided into four steps.\\
\textbf{Step 1. Construct a polynomial verifying the compatibility condition.}
We seek a polynomial condition for the selected squared distances in \eqref{eqn:sq-dist} to be realized by a single point. If $a'\in\R^d$ satisfies
\begin{equation*}
|a'-b_j|^2=r_j^\alpha(x),
\quad j=0,\ldots,d,
\end{equation*}
then subtracting the equation at $b_0$ from the others gives
\begin{equation}\label{eqn:syst-generic}
2(b_j-b_0)\cdot a' = |b_j|^2-|b_0|^2+r_0^\alpha(x)-r_j^\alpha(x), \quad j=1,\ldots,d.
\end{equation}
The coefficient matrix of the linear system \eqref{eqn:syst-generic} is $M=2[b_1-b_0,\ldots,b_d-b_0]^\top$, which is invertible by affine independence of $b_0,\ldots,b_d$. Thus, \eqref{eqn:syst-generic} determines a unique candidate location, denoted by
\begin{equation*}
T_\alpha(x)
=
M^{-1}
\begin{pmatrix}
|b_1|^2-|b_0|^2+r_0^\alpha(x)-r_1^\alpha(x)\\
\vdots\\
|b_d|^2-|b_0|^2+r_0^\alpha(x)-r_d^\alpha(x)
\end{pmatrix}.
\end{equation*}
The candidate $T_\alpha(x)$ realizes all the selected squared distances if and only if it also satisfies the equation at $b_0$. Accordingly, we define the polynomial
\begin{equation*}
P_\alpha(x)=|T_\alpha(x)-b_0|^2-r_0^\alpha(x).
\end{equation*}
Thus, the selected squared distances are realizable by a single point if and only if $P_\alpha(x)=0$, in which case that point is $T_\alpha(x)$. Next we prove that $P_\alpha$ is not identically zero
whenever the indices in $\alpha$ are not all equal. We choose $k\in\{i_1,\ldots,i_d\}$ with $k\neq i_0$, and let $e\in\R^d$ record the positions at which $k$ occurs:
\begin{equation*}
e_j=
\begin{cases}
1,& i_j=k,\\
0,& i_j\neq k,
\end{cases}
\quad j=1,\ldots,d.
\end{equation*}
Then $e\neq0$. Fix all source positions except $x_k$. Since $k\neq i_0$, the quantity $r_0^\alpha(x)$ is independent of $x_k$. For each $j\geq1$ with $i_j=k$, we have $r_j^\alpha(x)=|x_k|^2-2b_j\cdot x_k+|b_j|^2$,
whereas the remaining $r_j^\alpha(x)$ are independent of $x_k$. Substituting these expansions into the definition of $T_\alpha$ gives
\begin{equation*}
T_\alpha(x)-b_0=-|x_k|^2M^{-1}e+L(x_k),
\end{equation*}
where $L(x_k)$ is a vector-valued polynomial of degree at most one in $x_k$. It follows that
\begin{equation*}
P_\alpha(x)
=
|M^{-1}e|^2|x_k|^4
-2|x_k|^2(M^{-1}e)\cdot L(x_k)
+|L(x_k)|^2-r_0^\alpha(x).
\end{equation*}
Since $M$ is invertible and $e\neq0$, the coefficient $|M^{-1}e|^2$ is strictly positive. Thus the degree-four homogeneous part of $P_\alpha$ in $x_k$ is nonzero, proving that $P_\alpha$ is not identically zero.\\
\textbf{Step 2. Construct two exceptional sets.}
Now we exclude two types of exceptional configurations. First, to exclude a single point realizing distances selected from different sources, define
\begin{equation*}
E_n^{(1)}(B)
:=
\bigcup_{\substack{
\alpha=(i_0,\ldots,i_d)\in\{1,\ldots,n\}^{d+1}\\
i_0,\ldots,i_d\ \text{not all equal}
}}
\{x\in X_n:P_\alpha(x)=0\}.
\end{equation*}
Outside this set, if a point $a'\in\R^d$ satisfies $|a'-b_j|=|x_{i_j}-b_j|$ for every $j=0,\ldots,d$, then the indices $i_0,\ldots,i_d$ must all be equal. Second, to ensure that distinct sources have distinct distances from $b_0$, define
\begin{equation*}
E_n^{(2)}(B)
:=
\bigcup_{1\leq i<k\leq n}
\{x\in X_n:|x_i-b_0|^2-|x_k-b_0|^2=0\}.
\end{equation*}
The polynomials defining $E_n^{(1)}(B)$ are nonzero by the preceding argument. Those defining $E_n^{(2)}(B)$ are also nonzero, since $|x_i-b_0|^2-|x_k-b_0|^2$ has quadratic part $|x_i|^2-|x_k|^2$ for $i\neq k$. Each set is therefore a finite union of zero sets of nonzero real polynomials and hence has Lebesgue measure zero. Combining these two types of exceptions gives $E_n(B):=E_n^{(1)}(B)\cup E_n^{(2)}(B)$.
Thus $E_n(B)$ is a Lebesgue null set in $X_n$ and is independent of the positive strengths. For $n=1$, both unions are empty, so $E_1(B)=\varnothing$.\\
\textbf{Step 3. Prove unique identification outside the exceptional set.}
We claim that, for every $x=(x_1,\ldots,x_n)\in X_n\setminus E_n(B)$ and every choice of positive strengths, the source configuration is uniquely determined by its heat traces at $B$ among all finite source configurations with positive strengths. To prove this claim, fix $x=(x_1,\ldots,x_n)\in X_n\setminus E_n(B)$ and equip $A=\{x_1,\ldots,x_n\}$ with arbitrary positive strengths $q_{x_i}$. Suppose that another finite source configuration $(A',q')$, also with positive strengths, produces the same heat traces at $B$. By Theorem \ref{thm:equivalence}, the corresponding weighted distance measures agree at every point of $B$. We first show that $A'\subseteq A$. Fix any $a'\in A'$. At each measurement point $b_j$, the weighted distance measure of $(A',q')$ has positive mass at $|a'-b_j|$. Equality of the measures therefore provides an index $i_j\in\{1,\ldots,n\}$ such that
\begin{equation*}
|a'-b_j|^2=|x_{i_j}-b_j|^2,
\quad j=0,\ldots,d.
\end{equation*}
Let $\alpha=(i_0,\ldots,i_d)$. By the construction of $T_\alpha$ and $P_\alpha$, these equations imply $a'=T_\alpha(x)$ and $P_\alpha(x)=0$. Since $x\notin E_n^{(1)}(B)$, the indices $i_0,\ldots,i_d$ must all be equal. If their common value is $i$, then $a'$ and $x_i$ have the same distances from $b_0,\ldots,b_d$. Subtracting the squared-distance equation at $b_0$ from those at $b_1,\ldots,b_d$ yields
\begin{equation*}
M(a'-x_i)=0.
\end{equation*}
Since $M$ is invertible, we obtain $a'=x_i$. As $a'\in A'$ was arbitrary, this proves $A'\subseteq A$. It remains to show that every point of $A$ belongs to $A'$ and that the corresponding strengths agree. Since $x\notin E_n^{(2)}(B)$, the distances from $b_0$ to the points of $A$ are pairwise distinct. In view of $A'\subseteq A$, the only possible source in $A'$ contributing mass at $|x_i-b_0|$ is $x_i$ itself. Equality of the weighted distance measures at $b_0$, together with $q_{x_i}>0$, therefore implies
\begin{equation*}
x_i\in A'
\quad\text{and}\quad
q'_{x_i}=q_{x_i},
\quad i=1,\ldots,n.
\end{equation*}
Thus $A'=A$ and all strengths agree, proving the claim.\\
\textbf{Step 4. Prove sharpness.}
Let $B\subset\R^d$ be any fixed set of fewer than $d+1$ measurement points, and choose an affine hyperplane $H$ containing $B$. As in the proof of Corollary~\ref{cor:singleton}, every source location $a\notin H$ and its reflection across $H$ give two distinct unit-strength single-source configurations with identical traces at all points of $B$. Since $H$ has Lebesgue measure zero, uniqueness fails for almost every single-source location. This proves sharpness.
\end{proof}

\begin{remark}
The proof applies to every finite cardinality with the same fixed measurement set $B$. The exceptional sets $E_n(B)$ may depend on $B$ and $n$, but not on the positive strengths. Since exceptional configurations are allowed, this result is consistent with the sharp worst-case threshold in
Theorem \ref{thm:main}.
\end{remark}

\begin{remark}
No exceptional configurations need to be excluded when either
$d=1$ or $n=1$. For $d=1$, Theorem~\ref{thm:main-1d} shows that
two distinct points identify every finite configuration,
so $E_n(B)$ may be taken to be empty for every $n\geq1$.
For $n=1$, the proof above gives $E_1(B)=\varnothing$, in agreement
with Corollary~\ref{cor:singleton}.
\end{remark}

\section{Conclusion}\label{sec:conclusion}

In this work, we studied the unique identification of finitely many sources for the heat equation from time traces at finitely many spatial points. For a known real-valued radial Schwartz profile with nonzero integral, we showed that equality of heat traces on any nonempty open positive-time interval is equivalent to equality of the corresponding weighted distance measures. This equivalence reduces the inverse source problem to a geometric identification problem involving both the source locations and their positive strengths and provides the basis for sharp measurement bounds. For $d\geq 2$, any $dN+1$ measurement points in affine general position uniquely identify every configuration with at most $N$ source locations. This bound is sharp: some sets of $dN$ measurement points in affine general position fail to distinguish two distinct $N$-source configurations, even with unit strengths. In one dimension, two distinct measurement points suffice, independently of the number of sources. We also showed that $d+1$ measurements suffice for clustered designs under a reference-distance separation condition, for adaptive sampling of arbitrary finite configurations, and for almost every location configuration under a fixed affinely independent design. In each case, uniqueness holds among all finite configurations with positive strengths. These results clarify the role of unlabeled distance matching and show how the required number of spatial measurements depends on the source configurations and the measurement design.

\bibliographystyle{abbrv}
\bibliography{point}

\end{document}